\documentclass[titlepage,11pt]{article}
\usepackage{latexsym}
\usepackage{amsfonts}
\usepackage{amsmath,enumerate}
\usepackage{tikz}
\usetikzlibrary{graphs}
\usepackage{float}
\usepackage{amssymb}

\newcommand\blackslug{\hbox{\hskip 1pt \vrule width 4pt height 8pt depth 1.5pt
        \hskip 1pt}}
\newcommand\bbox{\hfill \quad \blackslug \bigbreak}

\def\ll{,\ldots,}

\title{Short directed cycles in bipartite digraphs}
\author{Paul Seymour\thanks{Supported by ONR grant N00014-14-1-0084 and NSF
    grant DMS-1265563.} and
  Sophie Spirkl\thanks{Current address: Rutgers University, Piscataway Township, NJ 08854.}\\
Princeton University, Princeton, NJ 08544}

\date{August 7, 2018; revised \today}

\newtheorem{thm}{}[section]

\newcommand{\Proof}{\noindent{\bf Proof.}\ \ }

\begin{document}
\maketitle
\begin{abstract}
The Caccetta-H\"{a}ggkvist conjecture implies that for every integer $k\ge 1$, if $G$ is a bipartite 
digraph, with $n$ vertices in each part, and every vertex
has out-degree more than $n/(k+1)$, then $G$ has a directed cycle of length at most $2k$. If true this is best 
possible, and we prove this for $k = 1,2,3,4,6$ and all $k\ge 224{,}539$.

More generally, we conjecture that for every integer $k\ge 1$, and every pair of reals $\alpha, \beta> 0$ with $k\alpha +\beta>1$,
if $G$ is a bipartite digraph with bipartition $(A,B)$, where every vertex in $A$ has out-degree at least $\beta|B|$,
and every vertex in $B$ has out-degree at least $\alpha|A|$, then $G$ has a directed cycle of length at most $2k$. This implies
the Caccetta-H\"{a}ggkvist conjecture (set $\beta>0$ and very small), and again is best possible for infinitely many pairs $(\alpha,\beta)$. 
We prove this for 
$k = 1,2$, and prove a weaker statement (that $\alpha+\beta>2/(k+1)$ suffices) for $k=3,4$.
\end{abstract}

\section{Introduction}

The Caccetta-H\"{a}ggkvist conjecture~\cite{ch}
states:
\begin{thm}\label{CHconj}
{\bf Conjecture (\cite{ch}): } For every integer $k\ge 1$, and all $n>0$, 
every $n$-vertex digraph in which every vertex has out-degree at least $n/k$
has girth at most $k$.
\end{thm}
(Digraphs in this paper are finite, and without loops or parallel edges; thus, there may be an edge $uv$ and another edge $vu$,
but not two edges $uv$. A digraph has {\em girth at most $k$} if it has a directed cycle of length at most $k$.)
This is vacuous for $k = 1$ and trivial for $k=2$; but for $k\ge 3$ it remains open, and indeed for $k=3$ it is one of the 
most well-known open questions in graph theory. There are numerous partial results (see~\cite{blair} for a survey).

In this paper we study an analogue of \ref{CHconj} for bipartite digraphs. 
(A digraph is {\em bipartite} if the graph underlying it is bipartite,
and a {\em bipartition} of a digraph means a bipartition of this graph.) If we take disjoint sets 
$V_1\ll V_{2k+2}$, each of cardinality $t$ for some fixed $t>0$, and make a digraph with vertex set
$V_1\cup\cdots\cup V_{2k+2}$, by making every vertex in $V_{i}$ adjacent from every vertex in $V_{i-1}$, for $1\le i\le 2k+2$
(where $V_{0}$ means $V_{2k+2}$), we obtain a digraph with a bipartition into parts both of cardinality $(k+1)t$ ($=n$ say),
in which every vertex has out-degree $n/(k+1)$, and with no directed cycle of length at most $2k$. This seems to be an extremal
example, and motivates the following conjecture, the topic of this paper:
\begin{thm}\label{mainconj}
{\bf Conjecture: }For every integer $k\ge 1$, if $G$ is a bipartite digraph, with $n>0$ vertices in each part, and every vertex
has out-degree more than $n/(k+1)$, then $G$ has girth at most $2k$.
\end{thm}

We show in section \ref{sec:CHimplies} that this is implied by the Caccetta-H\"{a}ggkvist conjecture \ref{CHconj}. 
The latter is not known to be true for any $k\ge 3$, but there have been many theorems showing that for certain values of $k$,
there exists $c$ such that every non-null digraph $G$ with minimum out-degree at least $c|V(G)|$ has girth at most $k$, 
where $c$ is a little bigger than, 
but close to, the conjectured $1/k$. We call these ``CH-approximations''. There is a curious phenomenon, that a good 
CH-approximation for some value of $k$ can be 
enough to imply \ref{mainconj} for that value of $k$. Because of this, we are able to prove \ref{mainconj} for a number
of values of $k$:
\begin{thm}\label{mainthm}
For $k = 1,2,3,4,6$, and all $k\ge 224{,}539$, if $G$ is a bipartite digraph, with $n>0$ vertices in each part, and every vertex
has out-degree more than $n/(k+1)$, then $G$ has girth at most $2k$.
\end{thm}
The proof for $k=3$ is given in section \ref{sec:6gons}, for $k=4$ in section \ref{sec:8gons}, for $k\ge 224{,}539$ in 
section \ref{sec:largek}, and for $k=6$ it is sketched
in section \ref{sec:12gons}.
There is more evidence for \ref{mainconj}; we will prove:
\begin{thm}\label{eulerian}
For all $k\ge 1$, if $G$ is a bipartite digraph, with $n>0$ vertices in each part,
and every vertex
has in-degree and out-degree exactly $\alpha n$, where $\alpha>1/(k+1)$, then $G$ has girth at most $2k$.
\end{thm}

Second, we examine a conjecture that is stronger than \ref{mainconj} (and, indeed, implies \ref{CHconj}):
\begin{thm}\label{unbalancedconj}
{\bf Conjecture:} For every integer $k\ge 1$, and every pair of reals $\alpha, \beta> 0$ with $k\alpha +\beta>1$,
if $G$ is a non-null bipartite digraph with bipartition $(A,B)$, where every vertex in $A$ has out-degree at least $\beta|B|$,
and every vertex in $B$ has out-degree at least $\alpha|A|$, then $G$ has girth at most $2k$.
\end{thm}
For each value of $k$, this (if true) is best possible for infinitely many values of $\alpha,\beta$, as we discuss later; 
for $\alpha=\beta$
it is equivalent to \ref{mainconj}, and for $\beta\rightarrow 0$ it implies \ref{CHconj}. Thus it would be optimistic to hope 
to prove it for $k=3$. It is easy for $k = 1$, and we will show it holds for $k=2$, that is:
\begin{thm}\label{2unbalancedthm}
For every pair of reals $\alpha, \beta> 0$ with $2\alpha +\beta>1$,
if $G$ is a non-null bipartite digraph with bipartition $(A,B)$, where every vertex in $A$ has out-degree at least $\beta|B|$,
and every vertex in $B$ has out-degree at least $\alpha|A|$, then $G$ has girth at most four.
\end{thm}
This theorem is the fundamental tool from which everything else in the paper is derived.

Third, we investigate variants of \ref{unbalancedconj} for small $k$. We prove a weakening, 
which still contains the $k=3$ and $k=4$ cases of \ref{mainthm}:
\begin{thm}\label{3unbalancedthm}
For $k=3,4$ and every pair of reals $\alpha, \beta> 0$ with $\alpha +\beta\ge 2/(k+1)$,
if $G$ is a bipartite digraph with bipartition $(A,B)$, where every vertex in $A$ has out-degree at least $\beta|B|$,
and every vertex in $B$ has out-degree at least $\alpha|A|$, then $G$ has girth at most $2k$.
\end{thm}

\section{Preliminaries}

Let us fill in the details of a few straightforward things from the previous section. First, to see that \ref{unbalancedconj}
implies \ref{CHconj}, suppose $H$ is a counterexample to \ref{CHconj}. Thus $H$ is a digraph with $n>0$ vertices,
and for some integer $k>0$, every vertex of $H$ has out-degree at least $n/k$, and $H$ has girth at least $k+1$.
Let $V(H) = \{h_1\ll h_n\}$.
Take two disjoint sets $A=\{a_1\ll a_n\}$ and $B=\{b_1\ll b_n\}$, and make a digraph $G$ with vertex set $A\cup B$ in which
$a_i$ is adjacent to $b_i$ for each $i$, and for every edge $h_ih_j$ of $H$, $b_i$ is adjacent to $a_j$ in $G$.
Then $G$ has girth at least $2k+2$. Moreover, every vertex in $B$ has at least $n/k$ out-neighbours in $A$, and every
vertex in $A$ has exactly one out-neighbour in $B$. Now let $\alpha = (1-\alpha)/k$ and $\beta = 1/(2n)$.
Then every vertex in $B$ has more than $\alpha n$ out-neighbours in $A$, and every vertex in $A$ has more than $\beta n$
out-neighbours in $B$; and so we have a counterexample to \ref{unbalancedconj}, for the same value of $k$.

Let $G$ be a digraph with a bipartition $(A,B)$; if $\alpha, \beta \ge 0$ are reals, we say $G$ is 
{\em $(\alpha,\beta)$-compliant} or {\em complies with $(\alpha,\beta)$}  (via the bipartition $(A,B)$)
if $G$ is non-null and every vertex in $A$ has at least $\beta|B|$ out-neighbours in $B$ and every vertex in $B$ has at 
least $\alpha |A|$ 
out-neighbours in $A$. Given $\alpha, \beta>0$, for which $k$ is it true that every $(\alpha,\beta)$-compliant digraph has girth
at most $2k$? (We denote this by $(\alpha,\beta)\rightarrow 2k$.) 
Note that by replacing every vertex in $A$ by a set of $|B|$ new vertices, and each vertex in $B$ by a set of $|A|$ new vertices,
and making each pair of new vertices adjacent if the corresponding pair of old vertices was adjacent, we may assume that
$|A|=|B|=n$ when it is convenient. 

Here is a useful example. 
Take three integers $k,s,t>0$, and let $n=k(s+t-1)+1$.
Take
two disjoint sets $A=\{a_1\ll a_n\}$ and $B=\{b_1\ll b_n\}$, and make a digraph $G$ with vertex set $A\cup B$ in which
for $1\le i\le n$ (reading subscripts modulo $n$):
\begin{itemize}
\item $a_i$ is adjacent to $b_j$ for $i\le j\le i+s-1$; and
\item $a_i$ is adjacent from $b_j$ for $i-t\le j<i$,
\end{itemize}
Then the digraph just constructed complies with 
$$\left(\frac{t}{k(s+t-1)+1}, \frac{s}{k(s+t-1)+1}\right),$$ 
and 
it is easy to see that it has girth at least $2n/(s+t-1)>2k$ (because starting from a vertex $b_i$
and moving along two edges of $G$ brings us to one of $b_{i+1}\ll b_{i+s+t-1}$). Consequently
$$\left(\frac{t}{k(s+t-1)+1}, \frac{s}{k(s+t-1)+1}\right)\not\rightarrow 2k.$$
Not all these statements are needed. 
For instance, if we set $s=2,t=3,k=4$, we find that $(3/17,2/17)\not\rightarrow 8$; and if we set $s=1,t=1,k=4$ we find that
$(1/5,1/5)\not\rightarrow 8$. The second fact implies the first because $(1/5,1/5)$ dominates $(3/17,2/17)$ in every coordinate, and in general we only
need to consider
the pairs given by the construction that are maximal under domination. These pairs 
turn out to be those with $s=1$ or $t=1$, that is, the pairs $(t/(kt+1), 1/(kt+1))$ and 
$(1/(kt+1), k/(kt+1))$ for $t\ge 1$, and henceforth we confine our attention to these.

Let us fix $k=2$. Say $(\alpha,\beta)$ is {\em good} if 
$(\alpha,\beta)\rightarrow 4$, and {\em bad} otherwise.
Which pairs are good?
The example shows that the pairs $(1/(2t+1), t/(2t+1))$ and $(t/(2t+1), 1/(2t+1))$ 
for $t=1,2,\ldots$ are bad,
and therefore so are all pairs they dominate.  (See figure 1.)
All these maximal bad pairs lie on one of the lines
$\alpha+2\beta=1$, $2\alpha+\beta=1$; and \ref{2unbalancedthm} says that all pairs outside the border formed 
by these two lines are good.
There are pairs that satisfy $\alpha+2\beta\le 1$ and
$2\alpha+\beta\le 1$ that are not dominated by pairs given by the example, and for most such pairs we do not know 
whether they are 
good or bad. With Alex Scott, we proved that all points on the dotted lines are good except for the bad points
of the example; and indeed, for every such good point, there is an open disc of good points including it. (We omit the proof.)

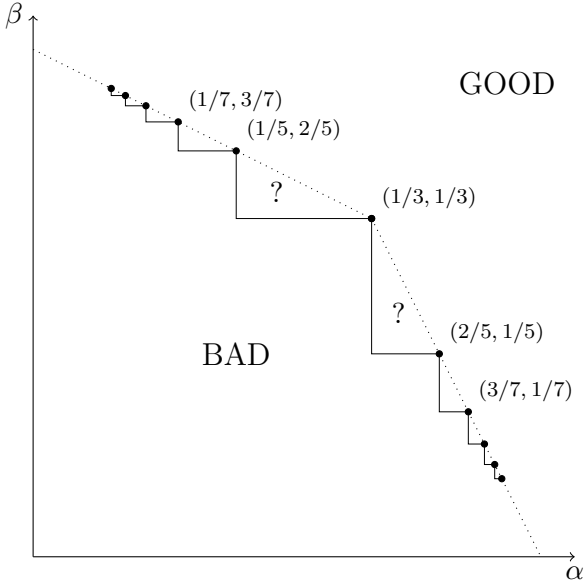
\begin{figure}[h]
\centering

\begin{tikzpicture}[scale=.9,auto=left]

\draw[->] (0,0) -- (8,0) node[anchor=north]{$\alpha$};
\draw[->] (0,0) -- (0,8) node[anchor=east] {$\beta$};
\foreach \x/\y in {5/5,3/6,2.1429/6.4286, 1.6667/6.6667, 1.36364/6.81818,1.15385/6.9231 }{
\draw[fill=black] (\x,\y) circle [radius=1.3pt];
\draw[fill=black] (\y,\x) circle [radius=1.3pt];}
\node[above right] at (5,5) {\begin{scriptsize}$(1/3,1/3)$\end{scriptsize}};
\node[above right] at (3,6) {\begin{scriptsize}$(1/5,2/5)$\end{scriptsize}};
\node[above right] at (2.1429,6.4286) {\begin{scriptsize}$(1/7,3/7)$\end{scriptsize}};
\node[above right] at (6,3) {\begin{scriptsize}$(2/5,1/5)$\end{scriptsize}};
\node[above right] at (6.4286,2.1429) {\begin{scriptsize}$(3/7,1/7)$\end{scriptsize}};
\draw[domain=0:5,smooth,variable=\x,dotted] plot ({\x},{7.5-\x/2});
\draw[domain=5:7.5,smooth,variable=\x,dotted] plot ({\x},{15-2*\x});
\draw[very thin] (5,5)-- (3,5) -- (3,6) -- (2.1429,6) -- (2.1429,6.4286) -- (1.6667,6.4286) -- (1.6667,6.6667) -- (1.36364,6.6667)
 -- (1.36364,6.81818) -- (1.15385,6.81818) -- (1.15385, 6.9231);
\draw[very thin] (5,5)-- (5,3) -- (6,3) -- (6,2.1429) -- (6.4286,2.1429) -- (6.4286,1.6667) -- (6.6667,1.6667) -- (6.6667,1.36364)
 -- (6.81818,1.36364) -- (6.81818,1.15385) -- (6.9231,1.15385);
\node at (3,3) {\large{BAD}};
\node at (7,7) { \large{GOOD}};
\node at (3.6,5.4) {$?$};
\node at (5.4,3.6) {$?$};

\end{tikzpicture}

\caption{The good, the bad, and the ugly.} \label{fig:1}
\end{figure}

For general $k$, much the same holds, except that for $k\ge 3$ we have not proved \ref{unbalancedconj}.
All the pairs given by the construction lie within the border given by the lines
$\alpha+k\beta=1$, $k\alpha+\beta=1$, and the maximal such pairs all lie on one of these lines; and we conjecture 
that all points $(\alpha,\beta)$ outside this border satisfy $(\alpha,\beta)\rightarrow 2k$.
This is one of the motivations for \ref{unbalancedconj} (the other is that it combines 
two interesting special cases, when $\alpha=\beta$ and when $\beta\rightarrow 0$).

Let us mention a difference between the Caccetta-H\"{a}ggkvist question and our bipartite question. 
The former asks whether minimum out-degree at least $c|V(G)|$ forces girth at most $k$ when $c=1/k$, 
and it is known that no smaller value of $c$ is big enough. For our bipartite question, we ask whether minimum out-degree
$c n$ forces girth at most $2k$, where $c>1/(k+1)$, and we know that $c= 1/(k+1)$ is not enough. (The point is, 
we need the strict inequality, while the Caccetta-H\"{a}ggkvist conjecture needs only the non-strict inequality.)

This causes a minor difficulty in 
terminology; for us, the important value is the maximal $c$ that is not ``big enough''. The pairs $(\alpha,\beta)$ we need to talk about are those such that $(\alpha, \beta) \not\rightarrow 2k$, and the pair $(\alpha, \beta)$ is maximal with that property, for the correct notion of maximality.

We could ask that every pair $(\alpha',\beta')$ that strictly dominates $(\alpha,\beta)$ satisfies $(\alpha',\beta')\rightarrow 2k$, where
``strictly dominates'' means $\alpha'\ge \alpha, \beta'\ge \beta$ and $\alpha'+\beta'>\alpha+\beta$; except this is not quite right.
When $\beta=0$, increasing only $\alpha$ will not guarantee any directed cycles at all, so in that case we need to say
$\beta'>\beta$; and we do not want to ask that both  $\alpha'> \alpha, \beta'> \beta$ hold, since that requires too much. Let us say
$(\alpha,\beta)\rightsquigarrow 2k$ if $(\alpha',\beta')\rightarrow 2k$ for every pair $\alpha',\beta'>0$ with
$\alpha'\ge \alpha, \beta'\ge \beta$ and $\alpha'+\beta'>\alpha+\beta$. Now we can write the statements from the first section more compactly; we can reformulate
\ref{CHconj} as ``for $k\ge 2$, $(1/k,0)\rightsquigarrow 2k$''; reformulate \ref{mainthm} as ``for $1\le k\le 4$,
$(1/(k+1),1/(k+1))\rightsquigarrow 2k$''; and reformulate \ref{unbalancedconj} as ``if $k\ge 1$ and $\alpha,\beta\ge 0$ with 
$k\alpha +\beta=1$, then
$(\alpha,\beta)\rightsquigarrow 2k$''.

The following implies the first ($k=1$) case of \ref{unbalancedconj}:
\begin{thm}\label{bipavdeg2}
Let $G$ be a digraph with a bipartition $(A,B)$, with more than $|A|\cdot|B|$ edges. Then $G$ has girth two.
\end{thm}
The proof is trivial and we omit it.

We need the following definitions: If $u,v\in V(G)$ where $G$ is a digraph, we say the {\em $G$-distance} from $u$ to $v$ is the number
of edges in the shortest directed path of $G$ from $u$ to $v$.
If $v$ is a vertex of a digraph $G$, $N_i(v)$ denotes the set of all vertices
$u$ such that the $G$-distance from $v$ to $u$ is exactly $i$, and
 $M_i(v)$ denotes the set of all vertices
$u$ such that the $G$-distance from $u$ to $v$ is exactly $i$.
Let $N^*_i(v)$ be the union of $N_j(v)$ over all $j$ with $1\le j\le i$ and with the same parity as $i$; thus $N_i^*(v)$
is a subset of one of $A,B$, if $(A,B)$ is a bipartition. We define $M_i^*(v)$ similarly.

\section{Previous results}

In our proofs, we will apply several previous results, each of which is a partial result towards proving the Caccetta-H\"aggkvist conjecture; and we have collected these results in this section.

We use these results in auxiliary digraphs; for example, in a bipartite digraph $G$ with bipartition $(A,B)$, we can pick $A' \subseteq A$ and construct a digraph $H$ with vertex set $A'$, and in which every vertex $v \in A'$ is adjacent in $H$ to $N_2(v)$ with respect to $G$. If $G$ has girth at least $2k$, then $H$ has girth at least $k$, and this allows us to deduce information about $G$ by applying an approximate version of the Caccetta-H\"aggkvist conjecture to $H$.

For this, we usually choose $A'$ such that every vertex in $A'$ has an out-neighbour in $B$ with all its out-neighbours in $A'$. Let $G$ be an $(\alpha,\beta)$-compliant digraph with $\alpha, \beta > 0$, and let $G$ have girth at least $2k$. If there is a vertex $v$ and $i \geq 3$ such that $N^*_i(v) \subseteq B$ and $|N_i(v)| < \beta |B|$, then we let $A' = N^*_{i-1}(v)$. Since every vertex in $N^*_{i-1}(v)$ has an out-neighbour in $N_{i-2}^*(v)$, it follows that every vertex in $H$ has $\alpha |A|$ out-neighbours in $H$. Therefore, if we know from an approximate version of the Caccetta-H\"aggkvist conjecture that every digraph $H$ of girth at least $k$ has a vertex of out-degree less than $|V(H)|/\delta$, then this construction shows that $\alpha |A| < |A'|/\delta$, and so $|N^*_{i-1}(v)| > \alpha\delta|A|$. We prove this formally in \ref{bigset}. 

The first result we use is the following, an approximation for the cases $k=3$ and $k=4$ of the Caccetta-H\"aggkvist conjecture.  
\begin{thm}[\cite{hladky, liang}] \label{CHbounds}
Let
$\delta_3=2.886$, and
let $\delta_4 = 3.4814$.
For $k=3,4$, every non-null digraph $G$ with minimum out-degree at least $|V(G)|/\delta_k$ has girth at most $k$.
\end{thm}
For $k=3$ this was proved by Hladk\'{y}, Kr\'al\textquoteright~and Norin~\cite{hladky}, and for $k=4$
it is due to Liang and Xu~\cite{liang}. (There is an improvement to $\delta_3= 2.9534$ claimed by Jan Volec in a 
lecture~\cite{volec}, but we do not use this.)

We also use the following result of Shen:
\begin{thm}[\cite{shen2}] \label{shen} 
Let $k\ge 1$ be an integer, and let $\delta=3k/4$;
then
every non-null digraph $H$ with minimum out-degree at least $|V(H)|/\delta$ has girth at most $k$.
\end{thm}
This is used in the proof of \ref{eulerian2}, an Eulerian version of \ref{mainconj}.

Finally, we need: 
\begin{thm}[\cite{shen}]\label{shenthm}
For $k>74$, if $H$ is a non-null digraph and every vertex has out-degree at least $|V(H)|/(k-74)$, then $H$ has girth
at most $k-1$.
\end{thm}
We use \ref{shenthm} in the proof of \ref{bigk} to say that for the first roughly $k-2r$ steps, the out-neighbourhood of a vertex grows steadily.

\section{Two results for all $k$}\label{sec:CHimplies}

In this section we prove \ref{eulerian}, and show that \ref{CHconj} implies \ref{mainconj}.

We need the following, which will be our main tool throughout the paper:

\begin{thm}\label{bigset}
Let $k\ge 1$ be an integer, and let $\delta>0$ be such that every non-null digraph $H$ with minimum out-degree at least $|V(H)|/\delta$ has girth at most $k$.
Let $\alpha,\beta> 0$, and let $G$ be an $(\alpha,\beta)$-compliant digraph
via a bipartition $(A,B)$, with girth more than $2k$.  Let $v\in V(G)$.
Then for each $i\ge 1$, if $N_i(v)\subseteq A$ then
either 
$|N_i(v)|\ge \alpha|A|$, or
$|N_{i-1}^*(v)|>\beta\delta |B|$; and if $N_i(v)\subseteq B$ then
either                
$|N_i(v)|\ge \beta|B|$, or
$|N_{i-1}^*(v)|>\alpha\delta |A|$.
\end{thm}
\Proof
From the symmetry we may assume that $N_i(v)\subseteq A$, and $|N_i(v)|< \alpha|A|$. Since $G$ is $(\alpha, \beta)$-compliant, it follows that $i\ge 2$.
Since $|N_i(v)|< \alpha|A|$, every vertex in $N_{i-1}^*(v)$ has an out-neighbour that does not belong to $N_i(v)$. In particular, every vertex in $N_{i-1}^*(v)$ has an out-neighbour in $N_{i-2}^*(v)$. Since $G$ has girth more than two, every vertex in $N_1(v)$ has all its out-neighbours distinct from $v$, and so $i\ge 3$. 
Let $H$ be the digraph with vertex set $N_{i-1}^*(v)$ in which for distinct $s,t\in N_{i-1}^*(v)$,
$t$ is adjacent from $s$ in $H$ if some vertex in $N_{i-2}^*(v)$ is adjacent in $G$ from $s$ and to $t$. 
Since every vertex in $N_{i-1}^*(v)$
has an out-neighbour in $N_{i-2}^*$, and every vertex in $N_{i-2}^*(v)$ has at least $\beta |B|$ out-neighbours in $N_{i-1}^*(v)$, it follows that
every vertex of $H$ has out-degree at least $\beta |B|$ in $H$. But $H$ has girth more than $k$, since $G$ has girth
more than $2k$; and so from the definition of $\delta$, since $V(H)\ne \emptyset$ it follows that in $H$, 
some vertex $u\in V(H)$ has
fewer than $|V(H)|/\delta$ out-neighbours in $V(H)$. Thus $|V(H)|>\beta \delta|B|$. This proves \ref{bigset}.~\bbox

Let us prove \ref{eulerian}, which we restate:
\begin{thm}\label{eulerian2}
Let $k\ge 1$ be an integer, let $\alpha>1/(k+1)$, and let $G$ be an $(\alpha,\alpha)$-compliant digraph via a
bipartition $(A,B)$, such that every vertex in $A$ has in-degree exactly $\alpha|B|$, and every vertex
in $B$ has in-degree exactly $\alpha|A|$. Then $G$ has girth at most $2k$.
\end{thm}
\Proof
We may assume that
$|A|=|B|=n$ say; and every vertex has in-degree exactly $\alpha n$ and out-degree exactly $\alpha n$.
There is therefore a perfect matching of edges each with tail in $A$ and head in $B$. 
\\
\\
(1) {\em $|N_i^*(v)| \ge ((i+1)/2)\alpha n$ for every vertex $v$ and for every odd integer $i\ge 1$ with $i\le 3k/2-1$.}
\\
\\
Suppose not and choose $i$ minimum such that this is false.
Thus $|N_i^*(v)|<((i+1)/2)\alpha n$. Note that $i \geq 3$. From the minimality of $i$, $|N_{i-2}^*(v)|\ge ((i-1)/2)\alpha n$, and so
$|N_{i}(v)|<\alpha n$. Let $\delta=3k/4$;
then by \ref{shen} and 
\ref{bigset}, $|N_{i-1}^*(v)|\ge \delta \alpha n$.
From the symmetry we may assume that $v\in A$, and so $N_i(v)\subseteq B$.
Since there is a perfect matching of edges each with tail in $A$
and head in $B$, and every vertex in $N_{i-1}^*(v)$ is matched into some vertex in $N_i^*(v)$, it follows that
$$((i+1)/2)\alpha n> |N_{i}^*(v)|\ge |N_{i-1}^*(v)|\ge \delta \alpha n,$$
and so $i > 2\delta-1=3k/2-1$, a contradiction. This proves (1).
\\
\\
(2) {\em If $k=2$ then $|N_3^*(v)| \ge 2\alpha n$ for every vertex $v$.}
\\
\\
Let $\delta=2$. Then every non-null digraph $H$ with minimum out-degree at least $|V(H)|/\delta$ has girth 
at most $2$, and now the claim follows as in (1). This proves (2).

\bigskip
If $k$ is odd, then since $k\le 3k/2-1$, it follows from (1) that $|N_k^*(v)|\ge ((k+1)/2)\alpha n>n/2$,
for every vertex $v$; and from \ref{bipavdeg2}, there are two vertices $u,v$ with $u\in N_k^*(v)$ and $v\in N_k^*(u)$,
and so $G$ has girth at most $2k$, and so \ref{eulerian} holds. So we assume that $k$ is even. Let $H$ be the digraph
with vertex set $V(G)$ in which for all $u\in A$ and $v\in B$, $u$ is adjacent to $v$ if $v\in N_{k-1}^*(u)$,
and $v$ is adjacent to $u$ if $u\in N_{k+1}^*(v)$. We claim that $H$ is $(((k+2)/2)\alpha, (k/2)\alpha)$-compliant:
if $k=2$ this follows from (2), and if $k\ge 4$ this follows from (1), since $k-1,k+1$ are odd, 
and both are at most $3k/2-1$.
Since $((k+2)/2)\alpha+(k/2)\alpha>1$,
the truth of \ref{unbalancedconj} for $k=1$ implies that $H$ has girth two, and so $G$ has girth at most $2k$.
This proves \ref{eulerian}.~\bbox

Next we show another of the claims of the introduction:

\noindent{\bf Proof of \ref{mainconj}, assuming \ref{CHconj}.\ \ }
Let $k\ge 1$ be an integer, let $\alpha>1/(k+1)$, and let $G$ be an $(\alpha,\alpha)$-compliant digraph via a bipartition $(A,B)$.
We must show that $G$ has girth at most $2k$. Suppose not. Since every vertex in $B$ has at least $\alpha|A|$ out-neighbours,
by averaging it follows that some vertex $v\in A$ has at least $\alpha|B|$ in-neighbours, that is, $|M_1(v)|\ge \alpha|B|$. 
Now the $k$ sets 
$N_1(v), N_3(v)\ll N_{2k-1}(v)$ are pairwise disjoint,
and all disjoint from $M_1(v)$ since $G$ has girth more than $2k$. Consequently not all of them have cardinality at least $\alpha|B|$,
since $|M_1(v)|\ge \alpha|B|$ and $|B| <(k+1)\alpha|B|$. By \ref{bigset}, taking $\delta=k$ (this 
satisfies the requirement of \ref{bigset} by the assumed truth of \ref{CHconj}) it follows that 
$|N_{2k-2}^*(v)|>\alpha k|A|.$
Since $N_{2k}(v)$ is disjoint from $N_{2k-2}^*(v)$, it follows that
$|N_{2k}(v)|\le |A|-k\alpha |A|<\alpha|A|$; and so by \ref{bigset} again,
$$|N_{2k-1}^*(v)|>k\alpha |B|>|B|-|M_1(v)|,$$
a contradiction. This proves \ref{mainconj}.~\bbox

We remark that the same proof shows the slightly stronger result that, assuming \ref{CHconj}, $(\alpha,\beta)\rightsquigarrow 2k$ 
for every pair $\alpha\ge \beta\ge 0$ with $\alpha+k\beta\ge 1$. (Note that this is very different 
from the conjecture \ref{unbalancedconj}.)

\section{Girth four}

In this section we prove \ref{2unbalancedthm}. 
We need the following lemma:

\begin{thm}\label{bigindeg}
For every pair of reals $\alpha, \beta> 0$,
if $G$ is an $(\alpha,\beta)$-compliant digraph with girth at least four, and with bipartition $(A,B)$, then for some vertex $v\in B$,
$|M_1(v)|+|M_3(v)|\ge (\alpha+\beta)|A|$.
\end{thm}
\Proof
Let $a=\lceil\alpha|A|\rceil$ and $b = \lceil\beta |B|\rceil$. Thus $b\ne 0$.
By deleting edges we may assume that every vertex in $A$ has out-degree
exactly $b$, and every vertex in $B$ has out-degree exactly $a$.
For each $v\in B$, let $w(v)= |M_1(v)|/(b|A|)$. Thus $\sum_{v\in B}w(v)=1$.
\\
\\
(1) {\em For each $u\in A$, $\sum_{v\in N_3(u)}w(v)\ge \alpha$.}
\\
\\
We write $N_i(u) = N_i$ for brevity. Thus $|N_1|=b$.
There are $a|N_1|=ab$ edges with tail in $N_1$, and they all have head in $N_2$. Since $G$ has girth at least four, it follows that there are at most
$b|N_2|-ab$ edges from $N_2$ to $N_1$. But there are $b|N_2|$ edges with tail
in $N_2$, so at least $ab$
edges from $N_2$ to $N_3$. Hence the sum of the in-degrees of the vertices in $N_3$ is at least $ab$; and so
$\sum_{v\in N_3(u)}w(v)\ge ab/(b|A|)\ge \alpha$. This proves (1).

\bigskip

Let us choose $v\in B$ at random, choosing each vertex $v$ with probability $w(v)$. The expectation of $w(v)$
is at least $1/|B|$
by the Cauchy-Schwarz inequality, and the expectation of $|M_1(v)|$
is $b|A|$ times the expectation of $w(v)$, and so at least $b|A|/|B|\ge \beta|A|$.
By (1), the expectation of $|M_3(v)|$ is at least $\alpha|A|$; and so the expectation of $|M_1(v)|+|M_3(v)|$
is at least $(\alpha+\beta)|A|$. Since $B\ne \emptyset$, some vertex $v\in B$ satisfies $|M_1(v)|+|M_3(v)|\ge (\alpha+\beta)|A|$.
This proves \ref{bigindeg}.~\bbox

We deduce \ref{2unbalancedthm}, which we restate:

\begin{thm}\label{2unbalancedthmagain}
$(\alpha,\beta)\rightsquigarrow 4$ for every pair of reals $\alpha, \beta\ge 0$ with $2\alpha +\beta\ge 1$.
\end{thm}
\Proof Let $\alpha'+\beta'>0$ with $2\alpha'+\beta'>1$ and $\beta'>0$; we must show that $(\alpha',\beta')\rightarrow 4$. Let $G$
be $(\alpha',\beta')$-compliant, via a bipartition $(A,B)$, and we assume for a contradiction that $G$ has girth at least six.
By \ref{bigindeg}, there exists $v\in B$
with $|M_1(v)|+|M_3(v)|\ge (\alpha'+\beta')|A|$. Since $|N_1(v)|\ge \alpha'|A|$, it follows that the sum of the cardinalities 
of the three sets $M_1(v), M_3(v)$
and $N_1(v)$ is at least $(2\alpha'+\beta')|A|>|A|$ (since $A\ne \emptyset$), and so some two of these three sets have nonempty
intersection, contradicting that $G$ has girth six. This proves \ref{2unbalancedthmagain}.~\bbox

We use $\mathbb{R}_+$ to denote the set of nonnegative real numbers. For an application later in the paper, we will need a more elaborate version of \ref{2unbalancedthm}, and to prove
that we need a lemma:

\begin{thm}\label{newineq}
Let $x,y,\beta,\gamma,\mu\in \mathbb{R}_+$ such that $x\le y\le 1$.
For $p,q,r\in \mathbb{R}_+$, define
$$f(p,q,r) = x(p-\gamma)^2 + (y-x)q^2 + (1-y)r^2.$$
Let $p,q,r\in \mathbb{R}_+$, with $px+q(y-x)+r(1-y) = \beta$ and $px+q(y-x) \ge \mu$.
\begin{enumerate}[(a)]
\item \label{it:one} If $\beta\le x\gamma$, then $f(p,q,r)\ge (\beta-x\gamma)^2/x$.
\item \label{it:two} If $\beta\ge x\gamma$, then $f(p,q,r)\ge (\beta-x\gamma)^2$.
\item \label{it:three} If $\beta\ge x\gamma$ and $y\beta +x(1-y)\gamma\le \mu$, then
$f(p,q,r)\ge (\mu- x\gamma)^2/y + (\beta-\mu)^2/(1-y).$
\end{enumerate}
(In each case, an expression with a zero denominator also has a zero numerator and should be taken to be zero.)
\end{thm}
\Proof

We first prove \eqref{it:one}. Since $x\gamma \geq \beta \geq px$, it follows that
$$f(p, q, r) \geq x(p-\gamma)^2 \geq (x\gamma-px)^2/x \geq (x\gamma - \beta)^2/x.$$

To prove \eqref{it:two}, we use the weighted inequality between quadratic mean and arithmetic mean that
$$xa^2 +(y-x)b^2 +(1-y)c^2 \geq (xa + (y-x)b + (1-y)c)^2$$
for all $a, b, c \in \mathbb{R}$ such that $xa + (y-x)b + (1-y)c \geq 0$. Since $0 \leq \beta - x\gamma = x(p-\gamma) + q(y-x)+r(1-y)$, setting $a = p-\gamma$, $b =q$, and $c = r$ implies that $f(p,q,r) \geq (\beta - x\gamma)^2.$

Finally, we prove \eqref{it:three}. The weighted inequality between quadratic mean and arithmetic mean implies that
$$xa^2 + (y-x)b^2 \geq (xa + (y-x)b)^2/y$$ for $a, b \in \mathbb{R}$ with $xa + (y-x)b \geq 0$. We let $a = p-\gamma$ and $b = q$. It follows that $xa + (y-x)b = \mu - x\gamma \geq y(\beta-x\gamma) \geq 0$, and hence $x(p-\gamma)^2 + (y-x)q^2 \geq  (\mu - x\gamma)^2/y$. Since $r \geq (\beta-\mu)/(1-y)$, it follows that $(1-y)r^2 \geq (\beta-\mu)^2/(1-y)$. Together, these inequalities imply that 
$$f(p,q,r) \geq (\mu - x\gamma)^2/y + (\beta-\mu)^2/(1-y).$$~\bbox

Next, we prove the \ref{appliedineq}, which in turn will be used to prove \ref{bellsandwhistles}. The aim of \ref{bellsandwhistles} is to take advantage of certain ``unbalanced'' parts of the digraph. For example, we will use \ref{appliedineq} with $b(v) |A| = |M_1(v)|$ and $a(v) = |N_{k-2}^*(v)|$ for some odd $k$. We will use the set $X$ below to capture vertices that can reach unusually many vertices in at most $k-2$ steps, and we will use $Y$ as a set with more incoming edges than the expected number $\beta |A||Y|$. 

\begin{thm}\label{appliedineq}
Let $B$ be a finite set, and let $x,y,\beta,\gamma,\lambda,\mu\in \mathbb{R}_+$ such that $x\le y\le 1$.
For each $v\in B$, let $a(v), b(v)\in \mathbb{R}_+$, such that:
\begin{itemize}
\item $\sum_{v\in B} b(v) = \beta |B|$;
\item for each $v\in B$, $a(v)\ge \lambda$;
\item there is a subset $X\subseteq B$ with $|X|= x|B|$ such that $a(v)\ge \gamma+\lambda$ for
all $v\in B\setminus X$.
\item there is a subset $Y\subseteq B$ with $|Y|= y|B|$, such that $\sum_{v\in Y} b(v)\ge \mu |B|$;
\item $\beta\ge x\gamma$ and $y\beta +x(1-y)\gamma\le \mu$.
\end{itemize}
Then $\sum_{v\in B}b(v)^2 + \sum_{v\in B}2a(v)b(v)\ge ((\mu- x\gamma)^2/y + (\beta-\mu)^2/(1-y) +2\beta (\lambda+\gamma)-x\gamma^2)|B| $.

(In each case, an expression with a zero denominator also has a zero numerator and should be taken to be zero.)
\end{thm}
\Proof
Let $g=\sum_{v\in B}b(v)^2 + \sum_{v\in B}2a(v)b(v)$.
We may assume that $a(v)=\lambda$ for each $v\in X$, and $a(v) = \lambda+\gamma$ for each $v\in B\setminus X$.
Consequently we may assume that $b(u)\ge b(v)$ for all $u\in X$ and $v\in B\setminus X$ (for if not, exchanging the
values of $b(u)$ and $b(v)$ reduces $g$). Consequently we may assume that $X\subseteq Y$ (since $y\ge x$).
Let 
\begin{eqnarray*}
\sum_{v\in X} b(v) &=& px|B|;\\
\sum_{v\in Y\setminus X} b(v) &=& q(y-x)|B|;\\
\sum_{v\in B\setminus Y} b(v) &=& r(1-y)|B|.
\end{eqnarray*} 
Thus $px+q(y-x) + r(1-y)=\beta$, and $px+q(y-x)\ge \mu$. 
By the Cauchy-Schwarz inequality, $\sum_{v\in X} b(v)^2\ge p^2x|B|$, and $\sum_{v\in Y\setminus X} b(v)^2\ge q^2(y-x)|B|$, and
$\sum_{v\in B\setminus Y} b(v)^2\ge r^2(1-y)|B|$; and so 
$$g/|B|\ge p^2x+ 2px\lambda + q^2(y-x) + 2q(y-x)(\lambda+\gamma) + r^2(1-y) + 2r(1-y)(\lambda+\gamma).$$
By \ref{newineq},
$$x(p-\gamma)^2 + (y-x)q^2 + (1-y)r^2\ge (\mu- x\gamma)^2/y + (\beta-\mu)^2/(1-y).$$
Thus
\begin{eqnarray*}
g/|B|&\ge& x(p^2+2p \lambda) + (y-x)(q^2+ 2q(\lambda+\gamma)) + (1-y)(r^2+2r(\lambda+\gamma)) \\
&\ge& (\mu- x\gamma)^2/y + (\beta-\mu)^2/(1-y) +(2px+2q(y-x)+2r(1-y)) (\lambda+\gamma)-x\gamma^2\\
&=&(\mu- x\gamma)^2/y + (\beta-\mu)^2/(1-y) +2\beta (\lambda+\gamma)-x\gamma^2.
\end{eqnarray*}
This proves \ref{appliedineq}.~\bbox

We now prove \ref{bellsandwhistles}. In the following, we will often take a bipartite digraph $G$ with bipartition $(A,B)$ and girth at least $2k$ for $k$ odd, and construct an auxiliary digraph $H$ with the same edges as in $G$ from $A$ to $B$, and in which every vertex $v$ in $B$ is adjacent to $N_{k-2}^*(v)$ with respect to $G$. Then $H$ has girth more than four, and we use methods similar to those in the proof of \ref{2unbalancedthm} to derive the desired inequality as one of the results of \ref{bellsandwhistles}. Furthermore, we generalize this approach to accommodate even numbers $k$ by introducing edge sets $R$ and $S$ as edges in $H$ from $v \in B$ to $N_{k-3}^*(v)$ in $G$ and $N_{k-1}^*(v)$ in $G$, respectively.

\begin{thm}\label{bellsandwhistles}
Let $x,y,\lambda,\beta,\gamma,\mu\ge 0$, with $x\le y\le 1$, and let $G$ be a digraph 
with a bipartition $(A,B)$.  Let $R,S \subseteq E(G)$, both consisting of edges of $G$ from $B$ to $A$. For $v\in B$, we denote
by $a_R(v)$ and $a_S(v)$ the number of edges in $R$ (respectively, $S$) with tail $v$; and let $a(v)= (a_R(v)+a_S(v))/(2|A|)$. Suppose the following hold.
\begin{itemize}
\item $\beta\ge x\gamma$ and $y\beta +x(1-y)\gamma\le \mu$.
\item Every vertex in $A$ has at least $\beta|B|$ out-neighbours in $B$.
\item $G$ has girth at least four and there is no directed cycle in $G$ of length four with an edge in $R$ and a 
different edge in $S$.
\item For each $v\in B$, $a(v)\ge \lambda$.
\item There is a subset $X\subseteq B$ with $|X|\le x|B|$ such that $a(v)\ge \gamma+\lambda$ for 
all $v\in B\setminus X$.
\item There is a subset $Y\subseteq B$ with $|Y|\le y|B|$, such that at least $\mu|A|\cdot |B|$ edges have head in $Y$.
\end{itemize}
Then
$$\frac{1}{y}(\mu-x\gamma)^2 + \frac{1}{1-y}(\beta-\mu)^2 + 2\beta(\lambda+\gamma)-x\gamma^2\le \beta.$$

(In each case, an expression with a zero denominator also has a zero numerator and should be taken to be zero.)
\end{thm}
\Proof We may assume that $\beta, x,y$ are rational; choose an integer $N>0$ such that $\beta N, xN, yN$ are 
all integers. Let $H$ be the digraph obtained from $G$ as follows; we replace each vertex $v$ of $G$ by a set 
$N_v$ of $N$
new vertices, and for each edge $uv$ of $G$ we add an edge from each member of $N_u$ to each member of $N_v$.
Let $R', S'$ be the sets of edges of $H$ that arise from edges in $R,S$ respectively. If the result holds for
$H$ then it holds for $G$; and so, replacing $G$ by $H$ if necessary, we may assume that $\beta|B|, x|B|, y|B|$ 
are all integers.

We may assume that every vertex in $A$ has out-degree exactly $\beta|B|$; and (by adding vertices to $X,Y$)
that $|X|=x|B|$ and $|Y|=y|B|$. 
We may assume that every edge from $B$ to $A$ belongs to $R\cup S$.
For each $v\in B$, let $b(v)|A|$ denote its in-degree. Thus $\sum_{v\in B} b(v) = \beta|B|$, since there
are $\beta|A|\cdot |B|$ edges from $A$ to $B$.

We define an {\em $R$-path} to be a directed path of $G$ such that all its edges with tail in $B$ belong to $R$; and an {\em $S$-path}
is defined similarly. 
For $u\in V(G)$ and $i\ge 1$, let $N_i^R(u)$ be the set of all vertices $v$ such that there is an $R$-path of length $i$
from $u$ to $v$, but no shorter $R$-path; and let $M_i^R(u)$ be the set of $v$ such that $u\in N_i^R(v)$. Define
$N_i^S(u), M_i^S(u)$ similarly.

For each edge $uv$ of $G$ with $u \in A$ and $v \in B$, the sets $M_1(v), M_3^R(v), N_1^S(v)$ are pairwise disjoint
subsets of $A$, and so the sum of their cardinalities is at most $|A|$. 
Summing over all $uv$, we deduce that
$$\sum_{v\in B}b(v)^2|A|^2 + \sum_{v\in B}b(v)a_S(v)|A| + \sum_{v\in B}b(v)|A|\cdot |M_3^R(v)| \le \beta |A|^2|B|,$$
 since there are
$\beta|A|\cdot |B|$ edges $uv$ from $A$ to $B$.
But $\sum_{v\in B}b(v)|A|\cdot |M_3^R(v)|$ is the sum over all $u\in A$, of the sum of the in-degrees of the vertices in $N_3^R(u)$;
and so at least the sum over $u\in A$, of the number of edges from $N_2^R(u)$ to $N_3^R(u)$. The latter is at least
the number of edges in $R$ from $N_1(u)$ to $N_2^R(u)$, since all vertices in $A$ have the same out-degree. The number of edges in $R$ from $N_1(u)$ to $N_2^R(u)$ is at least the sum
of $a_R(v)$ over all $v\in N_1(u)$.
So 
$$\sum_{v\in B}b(v)|A|\cdot |M_3^R(v)|\ge \sum_{v\in B} b(v)a_R(v)|A|.$$ 
The result therefore follows from \ref{appliedineq}.
This proves \ref{bellsandwhistles}.~\bbox

\section{Very large $k$}\label{sec:largek}

In this section we prove that \ref{mainconj} holds for all sufficiently large $k$.

\begin{thm}\label{bigk}
Let $r\ge 0$ be an integer such that for all $k\ge r+1$,
if $G$ is a non-null digraph and every vertex has out-degree at least $|V(G)|/(k-r)$, then $G$ has girth
at most $k-1$. Let $k$ be such that $k^2 +2(r+r^2)^2   + 2r^2    >  k(r^3/2+ 4r^2+4r)$ and $k\ge r(r+2)$.
(For instance with $r=74$, this requires $k> 224{,}538$.)
Then $(1/k,1/k)\rightsquigarrow 2k-2$.
\end{thm}
\Proof
Let $(A,B)$ be a bipartition of a digraph $G$, such that every vertex in $A$ has more than $|B|/k$ out-neighbours in $B$,
and every vertex in $B$ has more than $|A|/k$ out-neighbours in $A$. We must show that $G$ has girth at most $2k-2$.
Suppose not;
Inductively we assume the result for all digraphs with fewer vertices.
If $k$ is odd, for each $v\in B$ let $a(v)|A|=|N_{k-2}^*(v)|$. If $k$ is even let
$$a(v)|A|= (|N_{k-3}^*(v)|+ |N_{k-1}^*(v)|)/2= |N_{k-3}^*(v)|+|N_{k-1}(v)|/2.$$
Let $X$ be the set of all vertices $v\in B$ such that $a(v)\le (k-1)/(2k)$.
\\
\\
(1) {\em $X\ne \emptyset$.}
\\
\\
Let $H$ be the digraph with bipartition $A,B$, in which the edges from $A$ to $B$ are the same as in $G$,
and $u\in B$ is adjacent to $v\in A$ if the $G$-distance
from $u$ to $v$ is at most $k-1$. Let $S$ be the set of all edges of $H$ from $B$ to $A$. Let $R$ be the set of all edges $uv$ of
$H$ with $u\in B$ and $v\in A$
such that the $G$-distance from $u$ to $v$ is at most $k-2$. (So $S=R$ if $k$ is odd.) Thus, defining $a_R(v)$ and $a_S(v)$ as in
\ref{bellsandwhistles},
 we have $a_R(v)+a_S(v) = 2a(v)|A|$. Now the digraph $H$ has girth at least four, and has no directed 4-cycle
with an edge in $R$ and a different edge in $S$, since $G$ has girth at least $2k$. Let $\lambda$ be the minimum
of $a(v)$ over all $v\in B$. By \ref{bellsandwhistles}, taking
$x=0$, $y=1$, $\mu=\beta=1/k$, and $\gamma=0$, we have $1-y = 0 = \beta-\mu$, and therefore, 
$$\frac{1}{y}(\mu-x\gamma)^2 +  2\beta(\lambda+\gamma)-x\gamma^2\le \beta.$$
Consequently $1/k + 2\lambda\le 1.$ This proves (1).

\bigskip
Let $\delta = k-r$. From the choice of $r$, $\delta$ satisfies the condition of
\ref{bigset}.
\\
\\
(2) {\em There exists $Y\subseteq A$ with $|Y|\le |A|/2$ such that
at least $\delta|A|\cdot|B|/k^2$ edges have their head in $Y$. Moreover,
if $v\in X$ then
\begin{itemize}
\item  there are at least $\delta|B|/k$
vertices in $B$ with distance at most $k-2$ from $v$; and
\item $a(v)\ge (k-r-1)/(2k) $.
\end{itemize}}
\noindent Choose $v\in X$ (this is possible by (1)).
We claim that there exists $j$, odd, with $k-3\le j\le k-1$,
such that $|N_j^*(v)|\le (j+1)|A|/(2k)$. If $k$ is odd, this is true taking $j=k-2$, so
we assume that $k$ is even. Then
$$|N_{k-3}^*(v)|+|N_{k-1}^*(v)|=2a(v)|A|\le (k-1)|A|/k$$
and so either $|N_{k-3}^*(v)|\le (k-2)|A|/(2k)$ or $|N_{k-1}^*(v)|\le k|A|/(2k)$,
and again the claim holds, with $j=k-3$ or $k-1$.

Choose $g\ge 1$, odd,  minimal such that $|N_g^*(v)|\le (g+1)|A|/(2k)$. Thus $g\le k-1$, and
$|N_g^*(v)|\le (g+1)|A|/(2k)$, and $|N_{g-2}^*(v)| > (g-1)|A|/(2k)$,
so $|N_g(v)| \le |A|/k$.
By \ref{bigset}, $|N_{g-1}^*(v)|\ge \delta|B|/k= (1-r/k)|B|$. This proves the first bullet of (2).
In particular, $|N_g^*(v)|\le a(v)|A|\le |A|/2$;
all the edges with tail in $N_{g-1}^*(v)$ have head in $N_{g}^*(v)$, and since $|N_{g-1}^*(v)|\ge \delta|B|/k$,
this proves the first statement of (2) by taking $Y = N_g^*(v)$.

It remains to prove the second bullet of (2).
We claim that $|N_i(v)|\ge |A|/(2k)$ for all odd $i$ with $g\le i<k$. To see this, let $G'$ be the subdigraph
induced on
$N_{i-2}^*(v)\cup N_{i-1}^*(v)$. We may assume it is not $(\alpha,\alpha)$-compliant for any $\alpha>1/k$,
from the inductive hypothesis,
and since
$$|N_{i-2}^*(v)|\le |N_j^*(v)|< (j+1)|A|/(2k)\le |A|/2$$
it follows that some vertex in $N_{i-1}^*(v)$ has fewer than $|A|/(2k)$ out-neighbours in $N_{i-2}^*(v)$. Its other
out-neighbours all belong to $N_i(v)$, and so $|N_i(v)|\ge |A|/(2k)$, as claimed.

We claim that $a(v)\ge (k-r-1)/(2k)$. To see this there are two cases, depending on the parity of $k$. First, suppose
$k$ is odd.
Since
$$|N_{k-2}^*(v)|= a(v)|A|\le (k-1)|A|/(2k)<\delta |A|/k$$
(because $k> 2r-1$), it follows from \ref{bigset}
that $|N_i(v)|\ge |B|/k$ for all even $i$ with $1\le i\le k-1$. Since $N_{g-1}^*(v)$ is disjoint from the sets
$N_i(v)$ for all even $i$ with $g+1\le i\le k-1$, and there are $(k-g)/2$ such values of $i$, it follows that
$$(1-r/k)|B| + ((k-g)/2)|B|/k\le |B|,$$
that is, $k-g\le 2r$. In particular, $g\ge k-2r$.
So
$$a(v)= \frac{|N_{k-2}^*(v)|}{|A|}\ge \frac{|N_{g-2}^*(v)|}{|A|}+ \frac{(k-g)/2}{2k}\ge \frac{g-1}{2k}+\frac{k-g}{4k} \ge \frac{k-r-1}{2k}$$
as claimed, since $|N_{g-2}^*(v)|\ge (g-1)|A|/(2k)$, and $g\ge k-2r$.

Now we assume $k$ is even.
Since
$$|N_{k-3}^*(v)|< a(v)|A|\le (k-1)|A|/(2k)<\delta |A|/k$$
(because $k> 2r-1$), it follows from \ref{bigset}
that $|N_i(v)|\ge |B|/k$ for all even $i$ with $1\le i\le k-2$. Since $N_{g-1}^*(v)$ is disjoint from the sets
$N_i(v)$ for all even $i$ with $g+1\le i\le k-2$, it follows that
$$(1-r/k)|B| + ((k-1-g)/2)|B|/k\le |B|,$$
that is, $k-1-g\le 2r$. In particular, $g\ge k-1-2r$. As above, it follows that
$$\frac{|N_{k-3}^*(v)|}{|A|}\ge \frac{|N_{g-2}^*(v)|}{|A|}+ \frac{(k-g-1)/2}{2k}\ge \frac{g-1}{2k}+\frac{k-g-1}{4k} \ge \frac{k-r-2}{2k},$$
not quite what was promised. We can gain the missing $1/(2k)$ as follows. We may assume that $a(v)<(k-r-1)/(2k)$,
and since $|N_{k-3}^*(v)|/|A|\ge (k-r-2)/(2k)$, it follows that
$$\frac{|N_{k-1}^*(v)|}{|A|}\le \frac{k-r-1}{k} - \frac{k-r-2}{2k}= \frac{1}{2} -  \frac{r+1}{2k}<\delta/k;$$
and so by \ref{bigset}, $|N_k(v)|\ge |B|/k$. But then the argument above shows that
$$(1-r/k)|B| + ((k+1-g)/2)|B|/k\le |B|,$$ 
that is, $g\ge k+1-2r$; and then the estimate for $|N_{k-3}^*(v)|/|A|$
improves by the missing $1/(2k)$. This proves the claim that $a(v)\ge (k-r-1)/(2k)$,
and so proves (2).
\\
\\
(3) {\em $|X|\le 2r|B|/k$.}
\\
\\
Let $H$ be the digraph with vertex set $B$ in which $u$ is adjacent to $v$ if the $G$-distance from $u$ to $v$
is at most $k-1$. Then $H$ has girth at least three. Every vertex in $X$ has out-degree at least $\delta|B|/k$ in $H$;
and yet the number of edges of $H$ with tail in $X$ is at most $|X|^2/2 + |X|(|B|-|X|)$, since $H$ has no directed
cycles of length two. Hence $|X|\delta |B|/k\le |X|^2/2 + |X|(|B|-|X|)$, that is,
$|X|\le 2(1-\delta/k)|B|\le 2r|B|/k$. This proves (3).

\bigskip

From (1) and (2) and the symmetry between $A$ and $B$, there exists $Y\subseteq B$ with $|Y|\le |B|/2$ such that
at least $\delta|A|\cdot|B|/k^2$ edges have head in $Y$. By \ref{bellsandwhistles}, taking
$\lambda = (k-r-1)/(2k)$, $\beta=1/k$, $\gamma= r/(2k)$, $\mu=\delta/k^2$, $x=2r/k$ and $y=1/2$, we observe that
$\beta\ge x\gamma$  (since $k\ge r(r+2)$) and $y\beta +x(1-y)\gamma\le \mu$. By taking $R$ and $S$ as in the proof of (1), we deduce from \ref{bellsandwhistles} that
$$\frac{1}{y}(\mu-x\gamma)^2 + \frac{1}{1-y}(\beta-\mu)^2 + 2\beta(\lambda+\gamma)-x\gamma^2\le \beta,$$
that is,
$$2(\delta/k^2-(2r/k)r/(2k))^2 + 2(1/k-\delta/k^2)^2 + 2(1/k)((k-r-1)/(2k)+r/(2k))-(2r/k)r^2/(4k^2)\le 1/k,$$
which, recalling $\delta = k-r$, simplifies to
$$k^2 +2(r+r^2)^2   + 2r^2    \le  k(r^3/2+ 4r^2+4r),$$
a contradiction. This proves \ref{bigk}.~\bbox

\section{Girth six}\label{sec:6gons}

In this section we prove the case $k=3$ of \ref{3unbalancedthm}, that $(\alpha,\beta)\rightsquigarrow 6$
if $\alpha+\beta\ge 1/2$; and therefore the case $k=3$ of \ref{mainthm}. Equivalently, we must show:

\begin{thm}\label{3unbalanced}
$(\alpha,\beta)\rightarrow 6$ for all $\alpha,\beta>0$ with $\alpha+\beta>1/2$.
\end{thm}
\Proof Let $\alpha,\beta>0$ with $\alpha+\beta>1/2$, and let 
$G$ be $(\alpha,\beta)$-compliant, via a bipartition $(A,B)$; we will prove that $G$ has girth at most six, by induction on $|V(G)|$. By deleting edges, and by increasing $\alpha$ to $\lceil\alpha|A|\rceil / |A|$ and increasing $\beta$ to $\lceil\beta|B|\rceil / |B|$, we may assume that every vertex in $B$ has exactly $\alpha|A|$ out-neighbours in $A$, and every vertex in $A$ has exactly $\beta|B|$ out-neighbours in $B$.

We assume for a contradiction that $G$ has girth more than six.
Thus inductively, for all pairs $\alpha', \beta'>0$ with $\alpha'+\beta'>1/2$, no proper subdigraph of $G$ is 
$(\alpha',\beta')$-compliant. 
From the symmetry we may assume that $\alpha\ge \beta$. Let $\delta_3$ be as in \ref{CHbounds}.
\\
\\
(1) {\em There exists a real number $x$ with 
$$\max\left(\beta\delta_3,\frac{2\alpha\beta}{2\beta-\alpha}\right)\le x\le 1-\beta$$
such that some vertex $u\in A$ has in-degree at least $2\alpha x|B|$.}
\\
\\
By \ref{bigindeg}, there exists $v\in B$ with $|M_1(v)|+|M_3(v)|\ge (\alpha+\beta)|A|$; let $|N_2(v)|=x|B|$, and we will show 
that $x$ satisfies
the claim. Since $|N_1(v)| = \alpha|A|$, and the four sets $M_3(v), M_1(v), N_1(v), N_3(v)$ are pairwise disjoint subsets
of $A$, it follows that 
$$|N_3(v)|\le |A|-(2\alpha+\beta)|A|<\beta|A|\le \alpha|A|.$$ 
By \ref{bigset}, $|N_2(v)|\ge \beta\delta_3|B|$, and so $x\ge \beta\delta_3$. Since 
$$|N_1(v)\cup N_3(v)|\le |A|-(\alpha+\beta)|A|\le |A|/2< \alpha\delta_3|A|,$$
\ref{bigset} implies that $|N_4(v)|\ge \beta|B|$. Thus $|N_2(v)|\le |B|-\beta|B|$, and so $x\le 1-\beta$. Each vertex in $N_1(v)$
has at least $\beta|B|$ out-neighbours in $N_2(v)$; and
since $|N_3(v)|\le \beta|A|$, each vertex in $N_2(v)$ has at least $\alpha|A|-\beta|A| = ((\alpha-\beta)/\alpha)|N_1(v)|$
out-neighbours in $N_1(v)$. Consequently the subdigraph induced on $N_1(v)\cup N_2(v)$ complies with  
$$\left(\frac{\alpha-\beta}{\alpha},\frac{\beta|B|}{|N_2(v)|}\right),$$
and from the inductive hypothesis, 
$$\frac{\alpha-\beta}{\alpha}+\frac{\beta|B|}{|N_2(v)|}\le \frac{1}{2}.$$
Hence 
$$\frac{\beta|B|}{|N_2(v)|}\le \frac{1}{2}-\frac{\alpha-\beta}{\alpha}=\frac{\beta}{\alpha} -\frac{1}{2} = \frac{2\beta - \alpha}{\alpha},$$
and so
$(2\beta-\alpha)|N_2(v)|\ge 2\alpha\beta|B|$, and therefore $x\ge 2\alpha\beta/(2\beta-\alpha)$. 
Since there are at least $\alpha|A|\cdot|N_2(v)|$ edges with tail in $N_2(v)$, and they all have head in $N_1(v)\cup N_3(v)$,
and 
$$|N_1(v)\cup N_3(v)|\le |A|-(\alpha+\beta)|A|\le |A|/2,$$ 
some vertex in $N_1(v)\cup N_3(v)$ has in-degree at least
$\alpha|A||N_2(v)|/(|A|/2)=2\alpha x|B|$.
This proves~(1).
\\
\\
(2) {\em $\beta > 0.219$.}
\\
\\
From (1) it follows that $(2\beta-\alpha)(1-\beta)\ge 2\alpha\beta$, and so
since $\alpha>1/2-\beta$, it follows that
$(3\beta-1/2)(1-\beta)\ge (1-2\beta)\beta$
and so $\beta> 0.219$. This proves (2).

\bigskip
By (1), there exists  $u\in A$ with in-degree at least $2\alpha x|B|$.
\\
\\
(3) {\em $|N_1(u)\cup N_3(u)\cup N_5(u)|<\beta\delta_3|B|$.}
\\
\\
Since the sets $N_1(u)\cup N_3(u)\cup N_5(u)$ and $M_1(u)$ are disjoint, and $|M_1(u)|\ge 2\alpha x|B|$, it follows
that 
$$|N_1(u)\cup N_3(u)\cup N_5(u)|\le (1-2\alpha x)|B|\le (1-(1-2\beta) x)|B|.$$ 
Thus it suffices to show that 
$1-(1-2\beta) x <\beta\delta_3$, that is, $x>(1-\beta\delta_3)/(1-2\beta)$. 
Now from (1), $x\ge \beta\delta_3$, so we may assume that
$\beta\delta_3 \leq (1-\beta\delta_3)/(1-2\beta)$, 
and so $\beta<0.223$; and hence $0.219<\beta<0.223$ from (2). On the other hand, the second lower bound from (1) tells us
that $x\ge 2\alpha\beta/(2\beta-\alpha)$, that is, $x\ge (1-2\beta)\beta/(3\beta-1/2)$, and for $0.219<\beta<0.223$
$$(1-2\beta)\beta/(3\beta-1/2) \ge (1-\beta\delta_3)/(1-2\beta).$$
Consequently $x>(1-\beta\delta_3)/(1-2\beta)$
as required.  This proves (3).
\\
\\
(4) {\em $\beta> 0.242$, and $|N_6(u)|\ge \alpha|A|$, and $|N_1(u)|, |N_3(u)|\ge \beta |B|$.}
\\
\\
From (3) and \ref{bigset}, each of $|N_4(u)|, |N_6(u)|\ge \alpha|A|$. In particular, since $N_2(u), N_4(u), N_6(u)$
are pairwise disjoint, it follows that $|N_2(u)|\le (1-2\alpha)|A|\le \alpha \delta_3|A|$; and so by \ref{bigset},
$|N_3(u)|\ge \beta|B|$. Since also $|N_1(u)|\ge \beta|B|$, it follows from (3) that
$|N_5(u)|\le \beta\delta_3|B|-2\beta|B|< \beta|B|$ (since $\delta_3< 3$), and so by \ref{bigset},
$|N_2(u)\cup N_4(u)|>\alpha \delta_3|A|$. Since $|N_6(u)|\ge \alpha|A|$, it follows that $\alpha \delta_3|A| + \alpha|A|\le |A|$,
and so $\alpha < 1/(1+\delta_3)< 0.258$. So $\beta> 0.242$. This proves (4).

\bigskip

Since $\beta> 0.242$ and $\alpha+\beta\ge 1/2$ it follows that $\alpha\beta\ge 0.0624$, and so
$$|M_1(u)|\ge  2\alpha x|B|\ge 2\alpha \beta\delta_3|B|\ge 0.36|B|.$$
Let $A'=N_2(u)\cup N_4(u)$ and $B' = N_1(u)\cup N_3(u)$. Now $|A'|\le (1-\alpha)|A|$ since $|N_6(u)|\ge \alpha|A|$; and
so every vertex in $B'$ has at least $\alpha|A|\ge \alpha|A'|/(1-\alpha)\ge |A'|/3$ out-neighbours in $A'$. Let $|N_5(u)| = y|B|$;
since $|M_1(u)|\ge 0.36|B|$, it follows that $|B'|\le (0.64-y)|B|$. But every vertex in $A'$
has at least $\beta |B|$ out-neighbours in $B$, and they all lie in $B'$ except for at most $y|B|$ of them; and so the 
subdigraph induced on $A'\cup B'$ complies with
$$\left(\frac{1}{3}, \frac{\beta-y}{0.64-y}\right).$$
From the inductive hypothesis,
these two numbers sum to at most $1/2$; so $(\beta-y)/(0.64-y)\le 1/6$, that is, $5y\ge 6\beta-0.64$. Since $|M_1(u)|\ge 0.36|B|$, 
and $|N_1(u)|, |N_3(u)|\ge \beta |B|$, it follows that
$$0.36 + 2\beta + (6\beta-0.64)/5\le 1,$$
that is, $\beta\le 0.24$, contrary to (4). This proves \ref{3unbalanced}.~\bbox

\section{Girth eight}\label{sec:8gons}

Now we prove the remaining part of \ref{3unbalancedthm}, namely that $(\alpha,\beta)\rightsquigarrow 8$ for every pair
$\alpha,\beta\ge 0$ with $\alpha+\beta\ge 2/5$. The method is much the same as that for \ref{3unbalanced}; we
prove that some vertex $v\in A$ has big in-degree, and then work with the cardinalities of the various sets $N_i(v)$,
using $\ref{bigset}$ and induction on $|V(G)|$. One difficulty is that now we have to use $\delta_4$ rather than $\delta_3$,
and the value given for $\delta_4$ seems to be further (in some vague sense) from the conjectured $4$ than $\delta_3$ 
is from $3$.

We will prove:

\begin{thm}\label{4unbalanced}
For all $\alpha,\beta>0$ with $\alpha+\beta>2/5$, every $(\alpha,\beta)$-compliant digraph has girth at most eight.
\end{thm}
\Proof Let $G$ be a digraph that is $(\alpha,\beta)$-compliant via a bipartition $(A,B)$, where $\alpha,\beta>0$ 
and $\alpha+\beta>2/5$, and suppose that $G$ has girth at least ten. Inductively, we may assume
that for all $\alpha',\beta'>0$ with $\alpha'+\beta'>2/5$ and $\alpha'\ge \beta'$, no proper subdigraph of $G$ is
$(\alpha',\beta')$-compliant. From the symmetry we may assume that $\alpha\ge \beta$, and (by reducing $\alpha,\beta$)
that $\beta\le 1/5$.
\\
\\
(1) {\em Every vertex in $A$ has in-degree at most $(1-\beta\delta_4)|B|$.}
\\
\\
Let $x= 1-\beta\delta_4$, and suppose some vertex $v\in A$
has in-degree more than $x|B|$. Since $x+\beta\delta_4=1$, the union of $N_1(v), N_3(v), N_5(v)$ and $N_7(v)$
has cardinality less than $\beta|B|\delta_4$; and so by \ref{bigset}, each of $N_2(v), N_4(v), N_6(v), N_8(v)$
has cardinality at least $\alpha|A|$. In particular, the union of $N_2(v)$ and $N_4(v)$ has cardinality at most 
$|A|(1-2\alpha)<\alpha\delta_4|A|$, and so by \ref{bigset}, $N_1(v), N_3(v), N_5(v)$ all have cardinality at least $\beta|B|$.
The sets $M_1(v), N_1(v), N_3(v), N_5(v), N_7(v)$ are pairwise disjoint, and since
$x+4\beta >1$ and $|N_1(v)|\ge \beta |B|$, not all of $N_3(v), N_5(v), N_7(v)$ have cardinality at least $\beta|B|$.
By \ref{bigset}, $|N_2(v)\cup N_4(v)\cup N_6(v)|>\alpha \delta_4|A|$. 
Let $A'=N_2(v)\cup N_4(v)\cup N_6(v)$
and $B' = N_1(v)\cup N_3(v)\cup N_5(v)$. Thus $|A'|\le |A|(1-\alpha)$, since $|N_8(v)|\ge \alpha|A|$; and
$|B'|\le (1-x-y)|B|$ where $|N_7(v)|=y|B|$. The subdigraph induced on $A'\cup B'$ complies with 
$$\left(\frac{\alpha}{1-\alpha}, \frac{\beta-y}{1-x-y}\right),$$ 
and so
$$\frac{\alpha}{1-\alpha}+\frac{\beta-y}{1-x-y}\le 2/5$$
from the inductive hypothesis. But $x+y+3\beta\le 1$, and so
$\frac{\beta-y}{1-x-y} \geq \frac{4\beta-1+x}{3\beta}$ by letting $a = 3\beta - 1+x+y \leq 0$, $b = \beta - y$, $c = 1-x-y$, and noting that $\frac{a+b}{a+c} \leq \frac{b}{c}$ since $ca+bc \leq ab+bc$ (because $b \leq c$ and $a \leq 0$). This implies that
$$\frac{\alpha}{1-\alpha}+\frac{4\beta-1+x}{3\beta}\le 2/5.$$
Since $\alpha+\beta>2/5$, this implies 
$x< \beta/5 + 9/(3+5\beta) - 2$, a contradiction, since
$x= 1-\beta\delta_4> \beta/5 + 9/(3+5\beta) - 2$ for all $\beta$ with $0 <  \beta\le 1/5$. This proves (1).
\\
\\
(2) {\em Let $z = (2/5-\beta)(3/5+1/(1-\beta))$.
Then $|N_1(v)\cup N_3(v)|\ge z|A|$ for every vertex $v\in B$.}
\\
\\
Let $v\in B$, and let $|N_1(v)\cup N_3(v)|=w|A|$. We need to show that $w\ge z$, so we may assume that $w\le 2\alpha$. Hence
by \ref{bigset}, $|N_4(v)|\ge \beta|B|$, and so $|N_2(v)|\le (1-\beta)|B|$. The subdigraph induced on $N_1(v)\cup N_2(v)$
complies with
$$\left(\frac{2\alpha-w}{\alpha}, \frac{\beta}{1-\beta}\right),$$
and so
$$\frac{2\alpha-w}{\alpha}+ \frac{\beta}{1-\beta}\le 2/5,$$
that is,
$$w\ge \left(\frac{2}{5}-\beta\right)\left(\frac{3}{5}+\frac{1}{1-\beta}\right).$$
This proves (2).
\\
\\
(3) {\em $0.17<\beta< 0.19$.}
\\
\\
Let $H$ be the digraph with vertex set $A\cup B$, in which the edges from $A$ to $B$ are the same as in $G$, and
$v\in B$ is adjacent to $u\in A$ if the $G$-distance from $v$ to $u$ is at most three. Then $H$ has girth at least six.
Let $d$ be the minimum out-degree in $H$ of vertices in $B$; then by
\ref{2unbalancedthm}, $2d/|A|+\beta\le 1$, and so $d\le (1-\beta)|A|/2$; and so in $G$, some vertex $v\in B$
satisfies $|N_1(v)\cup N_3(v)| \le (1-\beta)|A|/2<2\alpha|A|$. By (2),
$$(2/5-\beta)(3/5+1/(1-\beta))\le (1-\beta)/2,$$
and so $\beta> 0.17$. Also, by \ref{bigset}, $|N_2(v)|\ge \beta\delta_4|B|$, and so there are at least $\alpha\beta|A|\delta_4|B|$
edges from $N_2(v)$ to $N_1(v)\cup N_3(v)$. Since $|N_1(v)\cup N_3(v)| \le (1-\beta)|A|/2$, (1) implies
that
$$\alpha\beta|A|\delta_4|B|\le ((1-\beta)|A|/2)(1-\beta\delta_4)|B|,$$
that is,
$$(4/5-2\beta)\beta\le (1-\beta)(1/\delta_4-\beta),$$
and so $\beta<0.19$. This proves (3).
\\
\\
(4) {\em Let $y=5\beta/(3+5\beta)+3\beta/5$.
Then $|N_2(v)\cup N_4(v)|\ge y|B|$ for every vertex $v\in B$.}
\\
\\
Let $|N_2(v)\cup N_4(v)|=x|B|$, and let $A'=N_1(v)\cup N_3(v)$, and $B'=N_2(v)$. Suppose that $x < y$.  
Since 
$|N_4(v)|\le (x-\beta)|B|$, every vertex in $A'$ has at least 
$$\beta|B|-|N_4(v)|=(\beta-x)|B|+|B'|\ge \frac{2\beta-x}{\beta} |B'|$$
out-neighbours in $B'$. 
By \ref{bigset}, since $x < y < \beta \delta_4$, it follows that $|N_5(v)|\ge \alpha|A|$, and so $|A'|\le (1-\alpha)|A|$. Consequently every vertex in $B'$ has at least 
$\alpha|A|\ge \frac{\alpha}{1-\alpha}|A'|$ out-neighbours in $A'$. Thus the subdigraph induced on $A'\cup B'$
complies with 
$$\left(\frac{\alpha}{1-\alpha},\frac{2\beta-x}{\beta}\right),$$ 
and so
$$\frac{\alpha}{1-\alpha}+\frac{2\beta-x}{\beta}\le 2/5,$$
that is,
$$ x\ge \frac{5\beta}{3+5\beta}+\frac{3\beta}{5}=y,$$
a contradiction. 
This proves (4).
\\
\\
(5) {\em Let 
$$x= \frac{\beta(\delta_4+1)-1/2}{\beta(\delta_4+1)-y}.$$
There are at least $x|B|$ vertices $v\in B$ with $|N_3(v)|\ge \alpha|A|$.}
\\
\\
Let $H$ be the digraph with vertex set $B$, in which for all distinct $u,v$, $v$ is adjacent from $u$ in $H$ if
the $G$-distance from $u$ to $v$ is at most four. Then $H$ has girth at least three, since $G$ has girth more than eight;
and so $|E(H)|\le |B|^2/2$. Thus 
$$\sum_{u\in B}|N_2(u)\cup N_4(u)|\le |B|^2/2.$$
Let $X$ be the set of all vertices $v\in B$ with $|N_3(v)|\ge \alpha|A|$. For each $v\in B\setminus X$,
\ref{bigset} implies that $|N_2(v)|\ge \beta\delta_4|A|$, and $|N_4(v)|\ge \beta|B|$; and so
$|N_2(u)\cup N_4(u)|\ge \beta(\delta_4+1)|B|$ for each $u\in B\setminus X$. By (4), $|N_2(u)\cup N_4(u)|\ge y|B|$
for each $u\in X$; and so 
$$(|B|-|X|)\beta(\delta_4+1)|B| + |X| (y|B|)\le |B|^2/2.$$
By taking $x = |X|/|B|$ and dividing by $|B|^2$, we obtain
$$x(\beta(\delta_4+1)-y )\ge \beta(\delta_4+1)-\frac12.$$
This proves (5).

\bigskip

Since $y=(5\beta/(3+5\beta)+3\beta/5)$, it follows from (3) that $y\ge 0.32$; and since
$$z = (2/5-\beta)(3/5+1/(1-\beta)),$$ 
it follows from (3) that $z\ge 0.38$.
Consequently 
$$x\ge  x'=\frac{\beta(\delta_4+1)-1/2}{\beta(\delta_4+1)-0.32}.$$
Let us apply \ref{bellsandwhistles} to the digraph $H$ of (3) above and setting $\lambda = 0.38, \gamma = 2\alpha-0.38, x=1-x', y=1, \mu = \beta$. By (2), (4) and (5), it follows that $\alpha \leq 1/2$ and 
$$(2\alpha-0.38)x'\left(2-\frac{(2\alpha-0.38)(1-x')}{\beta}\right) + 2(0.38) + \beta \le 1.$$
But this is impossible for $0.17<\beta< 0.19$, contrary to (3).
This proves \ref{4unbalanced}.~\bbox

\section{Girth twelve}\label{sec:12gons}

There is one other claim of the introduction that is not yet proved, that $(1/7,1/7)\rightsquigarrow 12$.
We do not give the proof in full because it contains no new ideas; here is a sketch.

First, show that every non-null digraph with minimum out-degree $|V(H)|/\delta$ has girth at most six, where
$\delta = 5.219$; this follows from the result of \ref{CHbounds} that $\delta_3=2.886$.
Now let $G$ be a digraph with a bipartition $(A,B)$, where every vertex in $A$ has more than $|B|/7$ out-neighbours in $B$
and vice versa. By \ref{2unbalancedthm}, there is a vertex $u\in A$ with $|N_1(u)\cup N_3(u)\cup N_5(u)|$
of cardinality at most $3|B|/7$. By \ref{bigset}, $|N_2(u)\cup N_4(u)|> \delta |A|/7$. Consequently there
are a disproportionate number of edges entering $N_1(u)\cup N_3(u)\cup N_5(v)$, and we can apply \ref{bellsandwhistles}
(taking $x=0$). We deduce that there is a vertex $u\in A$ with $|N_1(u)\cup N_3(u)\cup N_5(u)|\le .3993|B|$,
and consequently there is a vertex $v\in B$ of in-degree at least $.2667|A|$. Now $.2667> 1-\delta/7$,
and so by \ref{bigset}, the sets $N_2(v), N_4(v), N_6(v), N_8(v), N_{10}(v), N_{12}(v)$ all have cardinality
at least $|B|/7$. So $N_2(v)\cup N_4(v)\cup N_6(v)\cup N_8(v)$ has cardinality less than $\delta|B|/7$, and so
$N_3(v), N_5(v), N_7(v), N_9(v)$ all have cardinality at least $|A|/7$. Now let $A' $ be the union of
$N_3(v), N_5(v), N_7(v), N_9(v)$, and $B'$ the union of $N_2(v), N_4(v), N_6(v), N_8(v), N_{10}(v)$; every vertex
in $A'$ is adjacent to more than $|B'|/7$ vertices in $B'$ and vice versa, so we can apply induction to this subgraph.
This completes the sketch.

\section*{Acknowledgement}
We would like to thank Alex Scott for his help in classifying the good and bad points of the plane mentioned
in section 2. We thank Farid Bouya for pointing out a mistake in an earlier version of this paper. We are thankful to the anonymous referees for their helpful suggestions which included a simplified proof for \ref{newineq}.

\end{document}